\documentclass{ifacconf}

\usepackage{graphicx}      
\usepackage{natbib}        
\usepackage{amsmath}
\usepackage{amssymb}
\usepackage{mathtools}
\usepackage{xurl}

\usepackage[ruled,vlined,linesnumbered]{algorithm2e}
\DontPrintSemicolon
\SetKwInOut{Input}{Input}
\SetKwInOut{Output}{Output}
\SetKw{KwEnd}{end}
\usepackage{float}
\begin{document}
\begin{frontmatter}

\title{An Efficient Method for the Optimal Control of Microgrids Under Uncertainties using Local Reduction\thanksref{footnoteinfo}} 

\thanks[footnoteinfo]{This work was jointly funded by EPSRC and SLB.}

\author[First]{Edoardo Scaccia} 
\author[First,Second]{Eric C. Kerrigan} 
\author[Third]{Anna Sadowska}

\address[First]{Department of Electrical and Electronic Engineering, Imperial College London, SW7 2AZ London, UK (e-mail: e.scaccia24@imperial.ac.uk).}
\address[Second]{Department of Aeronautics, Imperial College London, SW7 2AZ London, UK (e-mail: e.kerrigan@imperial.ac.uk).}
\address[Third]{SLB Cambridge Research, CB3 0EL Cambridge, UK (e-mail: asadowska@slb.com)}

\begin{abstract}                
The problem of optimal sizing and power scheduling in microgrids subject to uncertainties is well known to the control community. Commonly, the optimal control problem is cast as a mixed-integer program to model the logical constraints arising in energy storage systems, and is then solved approximately using numerical methods such as the scenario approach. In this paper, we propose and compare two formulations of a robust microgrid sizing and power scheduling optimal control problem with logical constraints and uncertainties in the user's power demand, solar power generation, grid electricity prices and battery efficiencies. The first formulation uses binary variables and big-M constraints, leading to a mixed-integer linear program. The second formulation casts the problem as a continuous nonlinear program through an exact smooth reformulation of the logical constraints, consisting of additional modelling variables and non-convex constraints. We then propose a novel local reduction algorithm, extending an existing method, to solve both problems. The two formulations are compared by evaluating the solutions returned by local reduction using \(100{,}000\)-sample Monte Carlo simulations and achieve promising results, with both averaging feasibility rates above \(90\%\).
\end{abstract}

\begin{keyword}
Energy management systems, Energy storage systems, Energy communities.
\end{keyword}

\end{frontmatter}

\section{Introduction}
The continuous increase in the world's population contributes to the rise in global energy consumption. The operation of residential and commercial buildings accounts for 30\% of global final energy consumption \citep{IEA2023Buildings}. The global electricity demand is growing at an annual average rate of 3.4\% \citep{IEA2024Electricity}, and in order to cope with this hunger for electrical energy, the world's energy production is shifting towards renewable generation. This is expected to grow by 2.7 times by 2030, reaching 17,000 TWh \citep{IEA2024Renewables}. The main difficulty with renewable power is its uncertain and intermittent nature, making it challenging to generate renewable power on-demand. This issue brings the need for energy storage systems to store the excess generated energy for later deployment, which requires automatic control strategies.

\begin{figure}
\begin{center}
\includegraphics[width=8.4cm, trim=3.2cm 3.1cm 3.1cm 3.1cm,clip]{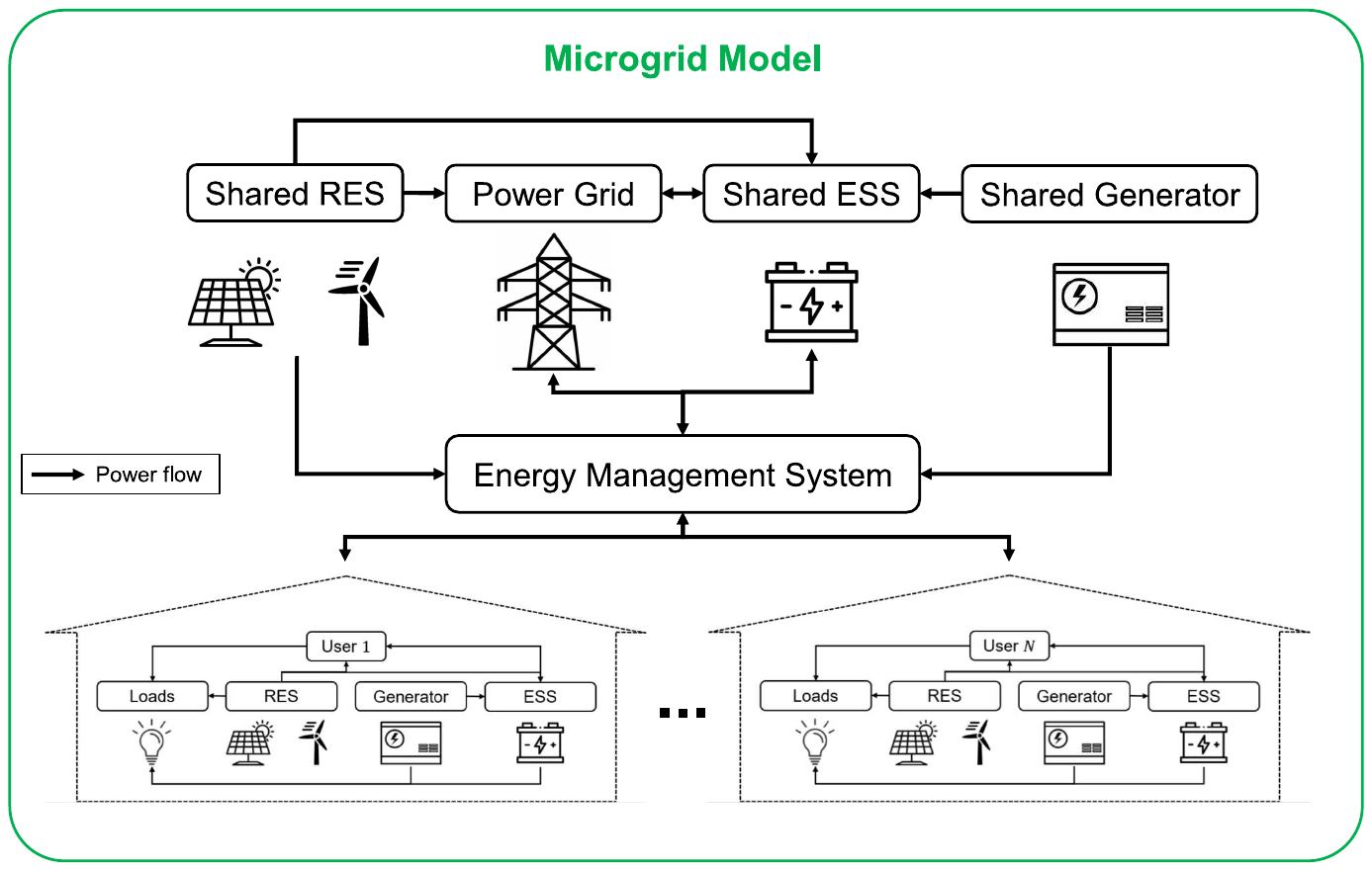}    
\caption{A general microgrid with multiple end users.} 
\label{fig:general_microgrid}
\end{center}
\end{figure}

Microgrids are localised electricity networks that use advanced control algorithms to generate, store and deploy electrical energy to meet the varying demand of end users while minimizing operational costs. A microgrid is comprised of: Energy Storage Systems (ESS), Renewable Energy Sources (RES) (e.g.\ wind and solar farms), non-renewable energy sources (e.g.\ dispatchable generators) and users (e.g.\ residential buildings, factories, offices, etc.). A schematic of a general microgrid is presented in Figure~\ref{fig:general_microgrid}. Microgrids are able to operate in a grid-connected or “islanded” mode. The connection to
the grid allows the system to trade power with the electricity market. On the other hand, when operating in islanded mode, microgrids are solely reliant on renewable and non-renewable energy sources to meet the users’ demand.

The optimal control of a microgrid is a challenging problem given the numerous uncertainties: renewable power generation, users’ power demand, grid's electricity prices, system parameters and output readings from sensors, system failure (e.g.\ outages, generators/wind turbines failures, etc.), ESS capacity degradation and maintenance schedules. Thus, an input sequence computed as the solution of the optimal control problem must be robust to every uncertainty realisation. This paper proposes a more efficient alternative to the scenario approach \citep{campi_scenario_2021} for solving this robust optimal control problem.

\section{State of the Art}
The problem of robust optimal control of microgrids is broken down into three main areas of focus: \textit{System Modelling}, \textit{Uncertainty Modelling} and \textit{Optimization Techniques}. Most microgrid optimal control formulations devote substantial modelling effort to the ESS, whose charge/discharge complementarity is commonly represented using binary variables and mixed-integer programming \citep{mousavizadeh_linear_2018}. RES, users' power demand and electricity prices are the main sources of uncertainty addressed in the literature, and are modelled employing three main approaches: min-max formulations \citep{zhang_robust_2013}, robust adaptable or ``multi-stage'' formulations \citep{wang_robust_2014} and probabilistic formulations \citep{mousavizadeh_linear_2018}. In min-max and robust adaptable problems, a robust optimal solution can be obtained by minimizing the worst-case scenario. However, in robust adaptable formulations, decision variables are split into first-stage (uncertainty-independent) and recourse decisions \citep{BenTal2009}. Lastly, in probabilistic formulations, uncertainties are modelled as random variables with associated probability density functions (PDFs). The most common PDFs employed in microgrid optimal control problems are the Gaussian/Normal, Weibull and Beta distributions. Probabilistic optimal control problems can minimize the expected value of the objective and can be characterised by \textit{chance} constraints \citep{Nemirovski2006}.

The most widely-used approach to solve linear and convex-concave quadratic min-max optimal control problems, is to use strong duality to substitute the dual of the inner max problem in the outer problem, resulting in standard min problems \citep{NocedalWright2006}. Min-max and probabilistic problems can also be approximately solved resorting to tractable approximations \citep{Nemirovski2006}. One common route is offered by the \textit{scenario approach} and Monte Carlo method \citep{yaghoubi_novel_2025}. Even though this approach imposes no restrictions on the distribution, a large number of samples are needed to obtain an accurate solution \citep{Nemirovski2006}. Lastly, linear robust adaptable two-stage problems are commonly solved by the ``Column-and-Constraint Generation'' method  \citep{zeng_solving_2013}, which requires strong assumptions on the uncertainty set.

\section{Contributions}
This paper considers a robust optimal sizing and scheduling problem for a grid-connected microgrid under uncertainty. The main contributions are fourfold. First, two formulations are proposed to encode the logical constraints on both battery and grid controls: an MILP formulation using binary variables and big-\(M\) constraints, and a continuous nonlinear programming (NLP) formulation using an exact smooth reformulation with auxiliary variables and non-convex constraints. Second, uncertainty is included directly in the ESS dynamics through the efficiencies, in addition to demand, solar power generation and electricity price uncertainty. Third, an extension of local reduction is used to correctly handle existence constraints via a logical reformulation of the lower-level maximization problem, leading to a nested local reduction scheme. Finally, the framework jointly optimizes the system size together with the robust battery scheduling decisions.

\section{Problem Formulation}
\label{sec:PF}
\begin{figure}
\begin{center}
\includegraphics[width=8.4cm, trim=4.9cm 3.6cm 7.8cm 1.1cm,clip]{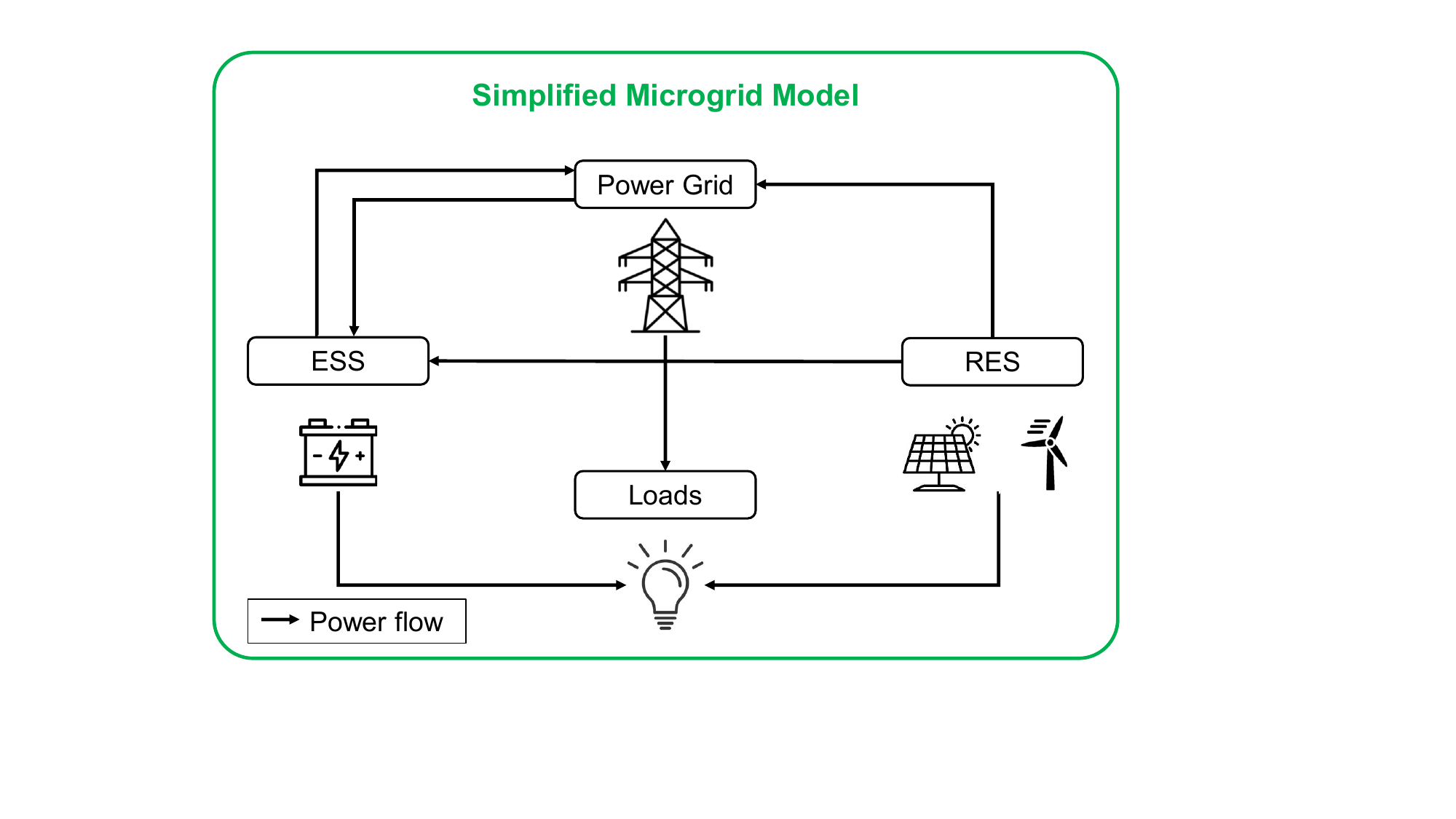}    
\caption{A schematic of a simplified microgrid with a single end user, ESS and RES.} 
\label{fig:Simple_Microgrid_New}
\end{center}
\end{figure}

A simplified microgrid model (Figure \ref{fig:Simple_Microgrid_New}), consisting of a single grid-connected user, ESS and RES, is considered. This formulation remains valid for more general models. The main model assumptions are introduced below, and the model parameters used in the numerical study are presented in Table \ref{tb:parameters}. Subsequently, the robust optimal control problem is formulated with the microgrid dynamics.

\subsection{Model Assumptions and Parameters}
\label{sec:assumptions}
\begin{enumerate}
  \item RES consists only of solar photovoltaic (PV).
  \item User's power demand $\tilde{P}_L$, PV power generation $\tilde{P}_R$, grid purchase price $\tilde{C}_{buy}$, grid sale price $\tilde{C}_{sell}$, ESS charging and discharging efficiencies $\tilde{\eta}_{\text{C}}$ and $\tilde{\eta}_{\text{D}}$ are uncertain.
  \item ESS dynamics are linear.
  \item ESS cannot charge and discharge simultaneously.
  \item User cannot buy and sell grid power simultaneously.
  \item Power losses, apart from ESS efficiency, are negligible.
  \item Line constraints are neglected.
  \item Model parameters and data are taken from \cite{vink_multiyear_2019} and modified when necessary. 
  \item ESS charging and discharging powers share the same limits.
  \item Grid import and export powers share the same limits. 
\end{enumerate}

\begin{table}[bt]
\begin{center}
\caption{Microgrid Parameters}
\label{tb:parameters}
\begin{tabular}{llcc}
\hline
\textbf{Parameter} & \textbf{Description} & \textbf{Value \& Units} \\
\hline
$SoC_0$             & Initial State of Charge    & $\left[0.2,1\right]$ \\
$SoC_{\min}$         & Minimum State of Charge    & 0.20 \\
$SoC_{\max}$         & Maximum State of Charge    & 1.00 \\
$P_{G,\max}$       & Maximum Grid Power          & - kW\\
$g$ & Grid Power Factor& $\left[0.2,1\right]$\\
$P_{ESS,\max}$       & Maximum ESS Power          & 90.0 kW\\
$E_{\text{ESS}}$     & ESS Base Capacity               & 326.0 kWh\\
$c_{\text{ESS}}$     & ESS Cost              & 200.00 £/kWh\\
$\bar{P}_{\text{PV}}$ & PV Rated Capacity & 90.84 kW \\
$c_{\text{PV}}$ & PV Cost & 1500.00 £/kW \\
$\tilde{C}_{buy,\text{UB}}$     & $\tilde{C}_{buy}$ Upper Bound (Vector)               &  $\text{£}/\text{kWh}$ \\
$\tilde{C}_{buy,\text{LB}}$     & $\tilde{C}_{buy}$ Lower Bound (Vector)               &  $\text{£}/\text{kWh}$ \\
$\tilde{C}_{sell,\text{UB}}$     & $\tilde{C}_{sell}$ Upper Bound (Vector)               &  $\text{£}/\text{kWh}$ \\
$\tilde{C}_{sell,\text{LB}}$     & $\tilde{C}_{sell}$ Lower Bound (Vector)               &  $\text{£}/\text{kWh}$ \\
$\tilde{P}_{L,\text{UB}}$     & $\tilde{P}_{L}$ Upper Bound (Vector)              & kW \\
$\tilde{P}_{L,\text{LB}}$     & $\tilde{P}_{L}$ Lower Bound (Vector)               & kW \\
$\tilde{P}_{R,\text{UB}}$     & $\tilde{P}_{R}$ Upper Bound (Vector)              & kW \\
$\tilde{P}_{R,\text{LB}}$     & $\tilde{P}_{R}$ Lower Bound (Vector)              & kW \\
$\tilde{\eta}_{\text{C,UB}}$   & $\tilde{\eta}_{\text{C}}$ Upper Bound (Vector)   & - \\
$\tilde{\eta}_{\text{C,LB}}$   & $\tilde{\eta}_{\text{C}}$ Lower Bound (Vector)    & - \\
$\tilde{\eta}_{\text{D,UB}}$  & $\tilde{\eta}_{\text{D}}$ Upper Bound (Vector) & - \\
$\tilde{\eta}_{\text{D,LB}}$  & $\tilde{\eta}_{\text{D}}$ Lower Bound (Vector) & - \\
\hline
\end{tabular}
\end{center}
\end{table}

\subsection{Microgrid Dynamics}
\label{sub:Dynamics}
The physics governing the simplified microgrid model is described by an ordinary differential equation (ODE) for the ESS energy level $E$ and by an algebraic equation
\begin{equation}
\label{eq:SoCODE}
\begin{aligned}
\frac{dE}{dt}= P_{c}(t)\cdot \tilde{\eta}_{\text{C}}(t)-P_{d}(t)\cdot \frac{1}{\tilde{\eta}_{\text{D}}(t)},
\end{aligned}
\end{equation}
\begin{equation}
\label{eq:Balance}
\begin{aligned}
    &P_{G}^{+}(t) - P_{G}^{-}(t) - P_{c}(t) + P_{d}(t) - \tilde{P}_L(t) + s \cdot \tilde{P}_R(t) = 0,
\end{aligned}
\end{equation}
where the grid power is split into its positive and negative components and the decision variable $s$ is a multiplicative factor applied to the base RES power to model additional solar panels. These differential-algebraic equations are complemented with the additional inequalities
\begin{align}
    &SoC_{\min} \cdot b \cdot E_{ESS} \le E(t) \le SoC_{\max} \cdot b \cdot E_{ESS},\label{eq:SoC_nounds}\\
    &0 \le P_{c}(t)  \le b \cdot P_{ESS,\max}, \label{eq:Pc_bounds}\\
    &0 \le P_{d}(t)  \le b \cdot P_{ESS,\max}, \label{eq:Pd_bounds}\\
    &0 \le P_{G}^+(t) \le g \cdot P_{G,\max},\label{eq:Positive_Grid_in}\\
    &0 \le P_{G}^-(t) \le g \cdot P_{G,\max}.\label{eq:Positive_Grid_out}
\end{align}
The decision variable $b$ represents a multiplicative factor applied to the ESS capacity and maximum powers to model additional battery packs. The parameter $g$ models the fraction of the allowable grid power that can be traded, where the base grid-power limit is set equal to the peak power demand over the considered horizon. Constraints derived from assumptions 4 and 5 are modelled using
\begin{equation}
\begin{aligned}
\label{eq:ESS_logical_statement}
&\textbf{IF} &P_{c}(t)>0 \\
&\textbf{THEN} &P_{d}(t) \le0, 
\end{aligned}
\end{equation}
\begin{equation}
\begin{aligned}
\label{eq:Grid_logical_statement}
&\textbf{IF} &P_{G}^{+}(t)>0 \\
&\textbf{THEN} &P_{G}^{-}(t) \le0. 
\end{aligned}
\end{equation}
The above discontinuous and non-smooth logical constraints need to be converted into equivalent continuous and smooth forms in order to be encoded into an optimization problem. Mixed-integer and non-convex smooth formulations are proposed. The former employs the big-M approach to reformulate \eqref{eq:ESS_logical_statement} and \eqref{eq:Grid_logical_statement} into
\begin{equation}
\begin{aligned}
\label{eq:ESS_bigM}
&0\le P_{c}(t) \le M_cx_1(t) \\
&0 \le P_{d}(t) \le M_d(1-x_1(t)),
\end{aligned}
\end{equation}
and
\begin{equation}
\begin{aligned}
\label{eq:Grid_bigM}
&0\le P_{G}^+(t) \le M_g^+x_2(t) \\
&0 \le P_{G}^-(t) \le M_g^-(1-x_2(t)).
\end{aligned}
\end{equation}
Here, $x_1(t)$ and $x_2(t)$ are binary variables and $M_c$, $M_d$, $M_g^+$ and $M_g^-$ are the big-M values. Since $P_{c}(t)$ and $P_{d}(t)$ have the same upper limit, $M_c = M_d=b_{\text{max}}\cdot P_{ESS,\text{max}}$, where $b_{\text{max}}$ is the maximum number of batteries allowed. Similarly, $P_{G}^+(t)$ and $P_{G}^-(t)$ share the same upper limit, thus, $M_g^+ = M_g^-=g\cdot P_{G,\text{max}}$. Note that while the constraints in \eqref{eq:Grid_bigM} replace \eqref{eq:Positive_Grid_in} and \eqref{eq:Positive_Grid_out}, this is not the case for \eqref{eq:ESS_bigM} as \eqref{eq:Pc_bounds} and \eqref{eq:Pd_bounds} are needed to avoid a bilinear big-M formulation due to $b$ being a decision variable and not a parameter as $g$. This formulation is referred to as the mixed-integer linear program (MILP) formulation. The equivalent smooth formulation is derived according to \cite{wehbeh_exact_2026} and results in
\begin{equation}
\label{eq:SmoothLogicConstraints}
\begin{aligned}
    &\lambda_1(t)\cdot P_{c}(t) +
    \left(1 - \lambda_1(t)\right) \cdot P_{d}(t) \le0 \\
    &\lambda_2(t)\cdot P_{G}^{+}(t)  +
    \left(1 - \lambda_2(t)\right) \cdot P_{G}^{-}(t) \le0
\end{aligned}
\end{equation}
where $0 \le\lambda_1(t),\lambda_2(t) \le1$. This formulation, referred to as NLP, requires constraints \eqref{eq:Pc_bounds}--\eqref{eq:Positive_Grid_out} to be retained.

\subsection{Robust Optimal Control Problem}
The ODE with expressions \eqref{eq:Balance}, \eqref{eq:SoC_nounds}--\eqref{eq:Positive_Grid_out} and \eqref{eq:ESS_bigM}--\eqref{eq:SmoothLogicConstraints} are now discretized. A one-day scheduling window going from $t_0=00\colon00$ to $t_f=24\colon00$ with 15-minute intervals $\Delta t$ is considered. The discretized grid is constituted by every discrete point of the window except the last one, such that $t_n=t_0+n\Delta t$ where $n\in \{0,\ldots,N=95\}$. Following this discretization, a variable evaluated at $t_n$ is denoted with the subscript $n$. A piecewise-constant parametrization is chosen for the disturbances and for the input trajectory, resulting in an open-loop control sequence. A piecewise-affine parametrization is chosen for the state trajectory. The ODE \eqref{eq:SoCODE} is discretized according to the \textit{Forward Euler Method}, resulting in an exact discretization. The integral of the \textit{Lagrange} term, the cost function of the optimal control problem, is discretized using the \textit{Left Riemann Sum} as a quadrature method. The objective is to minimize the initial investment and operational cost:
\begin{equation}
\label{eq:cost}
\begin{aligned}
&b \cdot c_{\text{ESS}}\cdot E_{\text{ESS}} + s \cdot c_{\text{PV}}\cdot \bar{P}_{\text{PV}} + \\ &\sum_{n = 0}^{N-1} (\tilde{C}_{buy,n} \cdot P_{G,n}^{+} -\tilde{C}_{sell,n} \cdot P_{G,n}^{-})\cdot\Delta t.
\end{aligned}
\end{equation}
Before stating the robust optimal control problem that we intend to solve, some notation is defined following \cite{wehbeh_semi-infinite_2024}. The vectors of physical states and uncertainties are defined as $z_{p,n} \coloneqq \left[E_n\right]^\top$ and $w_n \coloneqq [\tilde{P}_{L,n}\; \tilde{P}_{R,n}\; \tilde{C}_{buy,n}\; \tilde{C}_{sell,n}\; \tilde{\eta}_{\text{C},n}\; \tilde{\eta}_{\text{D},n}]^\top$ respectively. The vectors of modelling variables $z_{m,n}$ consist of two vectors $z_{g,n} \coloneqq \left[P_{G,n}^+\; P_{G,n}^-\;\right]^\top$ and $z_{l,n}$, representing the vector of grid variables and logic variables respectively. To differentiate the latter between the MILP and the NLP formulation, $z_{l,n}^{\text{MILP}} \coloneqq \left[x_{1,n}\; x_{2,n}\right]^\top$ and $z_{l,n}^{\text{NLP}} \coloneqq \left[\lambda_{1,n}\; \lambda_{2,n}\right]^\top$ are defined. Lastly, the vector of control variables is defined as $\theta_n \coloneqq [P_{c,n}\; P_{d,n}]^\top$. The vectors resulting from stacking the above over the entire horizon are denoted by $z_p$, $\theta$, $z_g$, $z_l$ and $w$, where $\theta$ also contains the variables $b$ and $s$. Their sizes are defined by $n_{z_{p}}$, $n_{\theta}$, $n_{w}$, $n_{z_g}$ and $n_{z_l}$. The controls and uncertainties belong to the sets $\Theta \subseteq \mathbb{R}^{n_{\theta}}$ and $\mathcal{W} \subseteq \mathbb{R}^{n_{w}}$ respectively. The vector $z_g$ belongs to the set $\mathcal{Z} \subseteq \mathbb{R}^{n_{z_g}}$. The vectors $z_l^{\text{MILP}}$ and $z_l^{\text{NLP}}$ belong to the sets $\mathcal{X} \subseteq \{0,1\}^{n_{z_l}}$ and $\Lambda \coloneqq \{\lambda\in\mathbb{R}^{n_{z_l}}: 0\le \lambda \le 1\}$. Lastly, the set $\mathcal{Z}_p(\theta, w)$ defines the admissible state trajectories as a function of the disturbances and can be derived by discretising the ODE with the initial condition $E_0=SoC_0 \cdot b \cdot E_{ESS}$ and the constraint \eqref{eq:SoC_nounds}. The uncertainty set is chosen as the polyhedron
\begin{equation}
\label{eq:OCP_uncertainty_set}
\mathcal{W} := \left\{ \, 
w \;\middle|\;
\begin{aligned}
   & w_{LB,n} \le w_{n} \le w_{UB,n} \\
   & \forall n \in \{0,\ldots,N-1\}
\end{aligned}
\right\}
\end{equation}
where $w_{LB}$ and $w_{UB}$ can be found in Table \ref{tb:parameters}. Note that any compact uncertainty set can be chosen, as the method presented remains valid. Moreover, the vector-valued functions defining the logic constraints in MILP and NLP form, and the objective are given by
\begin{equation}
\bar{g}_{\text{NLP}}(\theta,z_g,z_l^{\text{NLP}})
:=
\left[
\begin{array}{c}
\lambda_{1,n}\theta_{n,1}+(1-\lambda_{1,n})\theta_{n,2} \\
\lambda_{2,n}z_{g,n,1}+(1-\lambda_{2,n})z_{g,n,2}
\end{array}
\right]_{n=0}^{N-1},
\end{equation}
\begin{equation}
\bar{g}_{\text{MILP}}(\theta,z_g,z_l^{\text{MILP}})
:=
\left[
\begin{array}{c}
\theta_{n,1}-M_c x_{1,n} \\
\theta_{n,2}-M_d(1-x_{1,n}) \\
z_{g,n,1}-M_g^+ x_{2,n} \\
z_{g,n,2}-M_g^-(1-x_{2,n})
\end{array}
\right]_{n=0}^{N-1},
\end{equation}
\begin{equation}
\begin{aligned}
    \bar{f}(\theta, w, z_g)\coloneqq \;&b \cdot c_{\text{ESS}}\cdot E_{\text{ESS}} + s \cdot c_{\text{PV}}\cdot \bar{P}_{\text{PV}} + \\ &\sum_{n = 0}^{N-1}(w_{n,3} z_{g,n,1} - w_{n,4} z_{g,n,2})\cdot\Delta t.
\end{aligned}
\end{equation}
Noticing that $-g \cdot P_{G,\text{max}}\le z_{g,n,1} - z_{g,n,2} \le g \cdot P_{G,\text{max}}$, their difference is substituted using \eqref{eq:Balance}, and the function
\begin{equation}
\bar{g}_{\text{BAL}}(\theta,w)
:=
\left[
\begin{array}{c}
\theta_{n,2}-\theta_{n,1}+w_{n,1}-s w_{n,2}-gP_{G,\max} \\
-\theta_{n,2}+\theta_{n,1}-w_{n,1}+s w_{n,2}-gP_{G,\max}
\end{array}
\right]_{n=0}^{N-1}
\end{equation}
is defined. Lastly, to simplify the notation, the function
\begin{equation}
\bar{g}(\theta,w,z_g, z_l^i, \gamma) := \left[
\begin{aligned}
& \bar{g}_{i}(\theta,z_g,z_l^i)  \\ &\bar{g}_{\text{BAL}}(\theta,w) \\ &\bar{f}(\theta, w, z_g) - \gamma
\end{aligned}
\right]
\end{equation}
is introduced, where $i\in\{\text{NLP}, \text{MILP}\}$ and where $\gamma \in \Gamma \subseteq \mathbb{R}$ is an upper bound variable on the cost and the set $\Gamma$ can be chosen to be a bounded interval and thus a compact set. The robust optimal control problem is formulated as a min-max existence-constrained problem in epigraph form,
\begin{subequations}
\label{eq:OCP_master}
\begin{align}
    &\min_{\theta \in \Theta, \gamma \in \Gamma} \gamma \\
    & \; \text{s.t.} \quad \forall w \in \mathcal{W},\; \forall z_p \in \mathcal{Z}_p(\theta,w): \notag\\
    &\quad\quad\ \ \exists (z_g, z_l^i) \in \mathcal{Z} \times L:\; \bar{g}(\theta,w,z_g, z_l^i, \gamma) \le 0, \label{con:vec_val_exi}
\end{align}
\end{subequations}
where $L \coloneqq \Lambda$ if $i=\text{NLP}$ and $L \coloneqq \mathcal{X}$ if $i=\text{MILP}$. Problem~\eqref{eq:OCP_master} is an existence-constrained semi-infinite program (ECSIP) that can be reformulated into a standard semi-infinite program, satisfying all the assumptions in \citet{blankenship_infinitely_1976}, by following the method presented in \citet{wehbeh_semi-infinite_2024} and converting the existence quantifier into a min operator. However, this transformation cannot be achieved directly as $\bar{g}$ is a vector-valued function of dimension $n_{\bar{g}}$. In order to convert \eqref{con:vec_val_exi} into a scalar constraint, the maximum over all elements of $\bar{g}$ needs to be less than or equal to zero. After this, the $\exists$ quantifier can be replaced by a min over the existential variable and Problem~\eqref{eq:OCP_master} is equivalent to
\begin{subequations}
\label{eq:OCP_master_min_max}
\begin{align}
    &\min_{\theta \in \Theta, \gamma \in \Gamma} \gamma \\
    & \; \text{s.t.} \quad \forall w \in \mathcal{W},\; \forall z_p \in \mathcal{Z}_p(\theta,w): \notag\\
    &\quad\quad\ \ \min_{(z_g, z_l^i) \in \mathcal{Z} \times L} \max_{j=1,\ldots,n_{\bar{g}}} \bar{g}_j(\theta,w,z_g, z_l^i, \gamma) \le 0 \label{con:min_max_master}
\end{align}
\end{subequations}

\section{The Local Reduction Algorithm}
\label{sec:LR} 

The local reduction algorithm, first presented in \citet{blankenship_infinitely_1976}, consists of solving a sequence of minimization and maximization problems in order to iteratively build a finite uncertainty set that approximates the original uncountable set. Problem~\eqref{eq:OCP_master} is discretized using a finite number of scenarios $w_{p} \in \mathbb{W}_p$ where $p \in \{1,\ldots, k_1\}$ and is then solved, obtaining the solution $(\theta,\gamma)^*_{k_1}$. A new violating scenario is found by maximizing \eqref{con:min_max_master} over $w$ and $z_p$ at the current solution $(\theta,\gamma)^*_{k_1}$, solving
\begin{equation}
\label{eq:max_prblem_init}
\begin{aligned}
    \max_{w \in \mathcal{W}, z_p \in \mathcal{Z}_p(\theta^*,w)}\min_{(z_g, z_l^i) \in \mathcal{Z} \times L} \max_{j=1,\ldots,n_{\bar{g}}} \bar{g}_j(\theta^*,w,z_g, z_l^i, \gamma^*).
\end{aligned}
\end{equation}
Converting the above into its epigraph form and replacing the min with the $\forall$ quantifier, the problem results in
\begin{subequations}
\label{eq:max_epigraph}
\begin{align}
    &\max_{w \in \mathcal{W}, z_p \in \mathcal{Z}_p(\theta^*,w), \sigma \in \Sigma} \sigma \\
    & \; \text{s.t.} \quad \forall (z_g, z_l^i) \in \mathcal{Z} \times L: \notag \\
    &\quad\quad\ \ \sigma \le \max_{j=1,\ldots,n_{\bar{g}}}\bar{g}_j(\theta^*,w,z_g, z_l^i, \gamma^*), \label{con:sigma_max}
\end{align}
\end{subequations}
where $\sigma \in \Sigma$ is a lower bound variable on the cost and the set $\Sigma$ can be chosen as a compact interval in a similar fashion as per $\Gamma$. Note that \eqref{con:sigma_max} encodes a disjunction as at least one element of $\bar{g}$ must be greater than or equal to $\sigma$. Hence, employing the technique presented by \citet{wehbeh_exact_2026}, this disjunction is formulated as
\begin{subequations}
\label{eq:max_smooth_logic}
\begin{align}
    &\max_{w \in \mathcal{W}, z_p \in \mathcal{Z}_p(\theta^*,w), \sigma \in \Sigma} \sigma \\
    & \; \text{s.t.} \quad \forall (z_g, z_l^i) \in \mathcal{Z} \times L: \notag\\
    &\quad\quad\ \ \exists y \in \mathcal{Y}: \sum_{j=1}^{n_{\bar{g}}}y_j\cdot (\sigma - \bar{g}_j(\theta^*,w,z_g, z_l^i, \gamma^*)) \le 0 \label{con:sigma_max_smooth}
\end{align}
\end{subequations}
where $\mathcal{Y}\coloneqq \left\{y\in \mathbb{R}_{\ge0}^{n_{\bar{g}}}:\sum_{j=1}^{n_{\bar{g}}}y_j=1\right\}$. Alternatively, \eqref{con:sigma_max} can be reformulated by enforcing $y\in \{0,1\}^{n_{\bar{g}}}$ and by employing the big-M formulation
\begin{equation}
\label{big-M-subproblem}
\begin{aligned}
    \exists y \in \mathcal{Y}: \sigma \le \bar{g}_j(\theta^*,w,z_g, z_l^i, \gamma^*) + M_j(1-y_j)\; \forall j,
\end{aligned}
\end{equation}
where, for each added scenario $(z_g,z_l^i)$, the values $M_j$ are computed from $\max_{\sigma \in \Sigma, w\in \mathcal{W}}(\sigma - \bar{g}_j(\theta^*,w,z_g, z_l^i, \gamma^*))$. Thus, \eqref{eq:max_smooth_logic} is a nested ECSIP that can be solved again by means of the local reduction algorithm. The problem is discretized using a finite number of scenarios $(z_g,z_l^i)_q \in \mathbb{Z}_q$ where $q \in \{1,\ldots, k_2\}$ and it is then solved, obtaining the solution $(w,z_p,\sigma)^*_{k_2}$. Similarly, a violating scenario at the current solution is found by maximizing \eqref{con:sigma_max_smooth} over $z_g$ and $z_l^i$, obtaining the problem 
\begin{equation}
\label{eq:max_subproblem}
\begin{aligned}
    \max_{(z_g, z_l^i) \in \mathcal{Z} \times L}\min_{y \in \mathcal{Y}} \sum_{j=1}^{n_{\bar{g}}}y_j\cdot (\sigma - \bar{g}_j(\theta^*,w^*,z_g, z_l^i, \gamma^*)).
\end{aligned}
\end{equation}
Finally, converting \eqref{eq:max_subproblem} into its epigraph form, the min can be replaced by a $\forall$ quantifier over $y$. Since the set $\mathcal{Y}$ is a simplex, this is equivalent to enforcing the new epigraph variable to be less than or equal to all the elements of the function $\bar{g}$, leading to the problem
\begin{subequations}
\label{eq:max_sigma_bar}
\begin{align}
    &\max_{(z_g, z_l^i) \in \mathcal{Z} \times L, \bar{\sigma} \in \bar{\Sigma}} \bar{\sigma} \\
    & \; \text{s.t.} \quad \forall j \in \{1,\ldots,n_{\bar{g}}\}: \notag \\
    &\quad\quad\ \ \bar{\sigma}-(\sigma^* - \bar{g}_j(\theta^*,w^*,z_g, z_l^i, \gamma^*)) \le 0,
\end{align}
\end{subequations}
where $\bar{\sigma} \in \bar{\Sigma}$ is a lower bound variable on the cost and the set $\bar{\Sigma}$ can be chosen as a compact interval as done previously. For each local reduction loop, once the new violating scenario is found, it is added to the master problem which is then solved again. This procedure is repeated until either no more violating scenarios are found or a maximum number of scenarios are added to the master problem. The complete routine to solve Problem~\eqref{eq:OCP_master} is presented in Algorithm~\ref{alg:localred}.
\RestyleAlgo{ruled}
\begin{algorithm}[b]
\caption{Nested Local Reduction Algorithm}\label{alg:localred}
\Input{Maximum iterations: $k_{1,\text{max}}$ and $k_{2,\text{max}}$. Violation tolerances: $\text{tol}_1$ and $\text{tol}_2$. Initial uncertainty scenario guess $w_1$.}
\Output{Optimal solution $(\theta,\gamma)^*_{k_1}$ and the final set of scenarios $\mathbb{W}_{k_1}$.}
Initialize scenario sets: $\mathbb{W}_0 = \emptyset$ and $\mathbb{Z}_0 = \emptyset$.\\
Set $k_1\leftarrow 1$.\\
Initial uncertainty scenario set $\mathbb{W}_{k_1} \coloneqq \left\{w_{1}\right\}$.\\
\Repeat{$\lvert\mathbb{W}_{k_1}\rvert=\lvert\mathbb{W}_{k_{1}-1}\rvert$ \text{or} $k_1=k_{1,\text{max}}$}{
  \textit{Upper Level Local Reduction Loop}\\
  \textbf{Step 1.} Solve \eqref{eq:OCP_master} and obtain $(\theta,\gamma)^*_{k_1}$.\\
  \textbf{Step 2.} Compute initial scenario $(z_g,z_l^i)_1$. \\
  Set $k_2\leftarrow 1$.\\
  Initial uncertainty set $\mathbb{Z}_{k_2} \coloneqq \left\{(z_g,z_l^i)_1\right\}$.\\
\Repeat{$\lvert\mathbb{Z}_{k_2}\rvert=\lvert\mathbb{Z}_{k_2-1}\rvert$ \text{or} $k_2=k_{2,\text{max}}$}{
\textit{Lower Level Local Reduction Loop}\\
\textbf{Step 3.} Solve \eqref{eq:max_smooth_logic} using $(\theta,\gamma)^*_{k_1}$ and obtain $(w,z_p,\sigma)^*_{k_2}$.\\
\textbf{Step 4.} Solve \eqref{eq:max_sigma_bar} using $(\theta,\gamma)^*_{k_1}$, $(w,\sigma)^*_{k_2}$. \\
  \uIf{$\bar{\sigma}^*_{k_2} > \text{tol}_2$}{
  Add new scenario: $\mathbb{Z}_{k_{2}+1} \leftarrow \mathbb{Z}_{k_2} \cup \{(z_g,z_l^i)^*_{k_2}\}$. \\
  Set $k_2 \leftarrow k_2+1$.
  }
  \KwEnd \\
}
\uIf{$\sigma^*_{k_2} > \text{tol}_1$}{
  Add new scenario: $\mathbb{W}_{k_{1}+1} \leftarrow \mathbb{W}_{k_1} \cup \{w^*_{k_2}\}$. \\
  Set $k_1 \leftarrow k_1+1$.
  }
  \KwEnd \\
}
Return $(\theta,\gamma)^*_{k_1}$ and $\mathbb{W}_{k_1}$.
\end{algorithm}
The choice between the MILP and NLP formulations of the logical constraints changes the type of problems solved in Algorithm~\ref{alg:localred}. When employing the former and \eqref{big-M-subproblem} instead of \eqref{con:sigma_max_smooth}, a sequence of finite MILPs is solved in the local reduction routine. Assuming that the MILPs are solved to global optimality and that the local reduction terminates in a finite number of steps, the algorithm returns the global robust optimal solution. However, when employing the latter, Algorithm~\ref{alg:localred} solves a sequence of non-convex NLPs, thus local robust solutions are obtained.

\section{Results}
\label{sec:RES}
The numerical performance of Algorithm~\ref{alg:localred} is evaluated on Problem~\eqref{eq:OCP_master} under both the MILP and NLP formulations. The local reduction algorithm is implemented in the \texttt{Julia} programming language using the \texttt{JuMP} package and the \texttt{Ipopt} and \texttt{Gurobi} solvers. The implementation and data used to generate the numerical results are publicly available on GitHub.\footnote{\url{https://github.com/EdoScaccia/Local-Reduction-Algorithm-for-Optimal-Control-of-Microgrids}} The experiments are run on a Microsoft Surface Studio 2 laptop with a 13th Gen Intel\textsuperscript{\textregistered} Core\texttrademark{} i7 and 16~GB of RAM. The algorithm parameters are set as $k_{1,\text{max}}=10$, $k_{2,\text{max}}=10$, $\text{tol}_1=\text{tol}_2=10^{-6}$, and Problem~\eqref{eq:max_smooth_logic} in Step 3 is randomly initialised 5 times when solving the NLP formulation. The local reduction algorithm always converges before the last iteration when solving the MILP formulation, whereas it reaches the maximum number of scenarios for the NLP formulation.

For each of the NLP and MILP formulations, the grid power cap factor $g$ and the initial state of charge $SoC_0$ are varied independently over the range 20\%–100\%, resulting in 25 parameter combinations. The quality of the solutions is then assessed via Monte Carlo simulations with 100,000 uniformly randomly generated scenarios, and the statistics are presented in Table~\ref{tb:LRvsMC}. A scenario is deemed infeasible if it violates the $SoC$ bounds, grid power bounds, or logic constraints. Violations of the cost constraint are excluded from the feasibility count because they do not represent physical system infeasibilities; instead, they indicate that the computed upper bound $\gamma$ does not fully capture the realised worst-case cost. 

Both formulations achieve a similar feasibility score, but the MILP formulation is substantially faster than the NLP formulation. However, the NLP formulation leads to lower average constraint violations, with 0.099\% of total constraints violated compared with 0.56\%. When inspecting the constraint violations type breakdown, the most violated constraints are the cost and grid bounds for the MILP and NLP formulations respectively. Moreover, the MILP does not return any $SoC$ bounds violations compared to the NLP formulation that reports an average of 12.45\%. The NLP formulation results in a higher average $\gamma^*$, which explains the lack of cost violations and the lower maximum constraint violation of 1.54 that corresponds to a grid bound violation. The maximum violation resulting from the MILP formulation corresponds to a cost violation. Lastly, given that the MILPs in the local reduction algorithm are solved to global optimality, and thus no constraint violations are expected from the optimal solution, these violations might be attributed to tolerance values.

\begin{figure}
\begin{center}
\includegraphics[width=8.4cm, trim=0.2cm 0.2cm 0.3cm 1.2cm,clip]{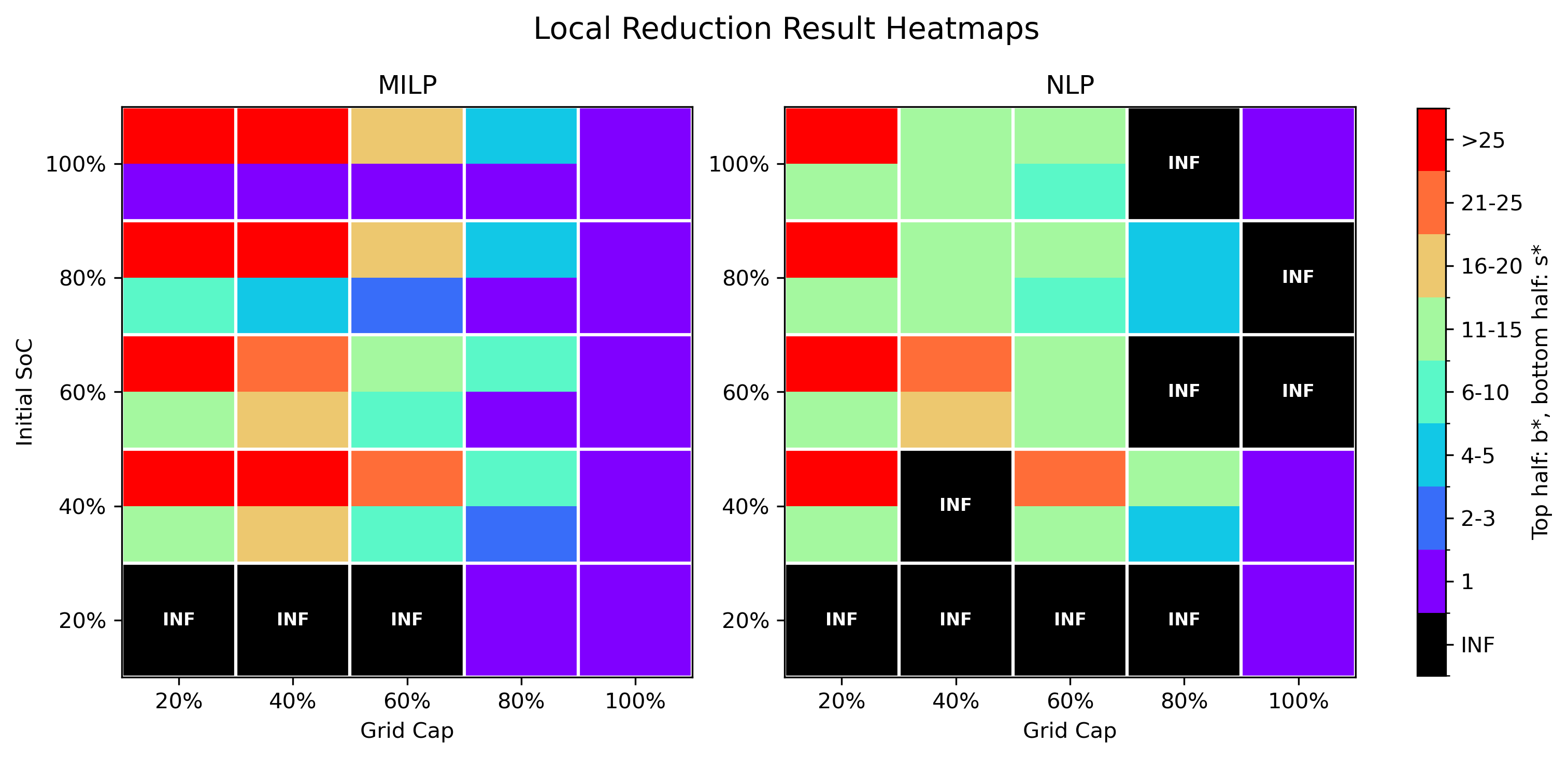}    
\caption{Optimal battery and PV sizing decisions as a function of the grid cap and initial \(SoC\).} 
\label{fig:LR_heatmaps}
\end{center}
\end{figure}

Figure~\ref{fig:LR_heatmaps} shows the optimal sizing decisions. The colour of the top half of each cell denotes the optimal number of battery packs $b^*$, while the bottom half reports the optimal number of PV units $s^*$; black cells indicate parameter combinations for which no feasible solution was returned. In the MILP formulation, the optimal decision variables $b^*$ and $s^*$ are integers, whereas in the NLP formulation they are continuous and are then mapped to integers by taking their ceiling values. The MILP formulation yields feasible solutions over a larger portion of the parameter grid, whereas the NLP formulation results in more infeasible cases. A notable difference is observed in the upper left part of the grid. The NLP formulation returns a more balanced number of battery packs and PV units compared to the MILP formulation, where the number of batteries is much higher compared to the number of PV units.

The difference in computational speed can be attributed to the structure of the optimization problems and the solvers used. In the MILP formulation, Gurobi, a state-of-the-art solver, solves linear problems. By contrast, the non-convex bilinear NLPs are solved using a general-purpose interior point solver. Significant speed-ups may be achievable through the development of a specialized solver tailored to the structure of Problem~\eqref{eq:max_smooth_logic}.

\begin{table}[tb]
\begin{center}
\caption{Monte Carlo simulations statistics of the MILP and NLP local reduction solutions.}\label{tb:LRvsMC}
\begin{tabular}{lcc}
& MILP & NLP \\\hline
Average Runtime (s) & 38.34 & 1410.15 \\
Maximum Runtime (s) & 282.244 & 2438.50\\ 
Average Feasible Scenarios (\%) & 90.91 & 90.65 \\
Average Constraints Violated (\%) & 0.56 & 0.099 \\
Average Violated SoC Bounds (\%) & 0.0 & 12.45 \\
Average Violated Logic (\%) & 0.0 & 0.0 \\ 
Average Violated Grid Bounds (\%) & 8.96 & 62.54 \\
Average Violated Cost (\%) & 72.86 & 0.0 \\
Average Maximum Violation (-) & 336.83 & 1.54 \\
Average $\gamma^*\cdot10^{-6}$ (£) & 2.09 & 2.64
\\\hline
\end{tabular}
\end{center}
\end{table}

\section{Conclusions}
This paper applied local reduction and an exact smooth reformulation of logical constraints to a robust microgrid sizing and scheduling problem. It is the first time that a microgrid optimal control problem under uncertainty is formulated as an NLP with embedded logic and solved via local reduction. The results show that both the proposed MILP and NLP formulations achieve average physical feasibility rates above \(90\%\) in \(100{,}000\)-sample Monte Carlo simulations. The method required at most \(10\) generated scenarios, whereas full enumeration of the box uncertainty set would require \(2^{6N}=2^{570}\) extreme scenarios. Assuming linear scaling, a deterministic equivalent with \(100{,}000\) scenarios would require roughly four orders of magnitude more memory and runtime than the corresponding \(10\)-scenario formulation. This illustrates the advantage of generating only relevant scenarios using local reduction.

Looking ahead, several avenues for future research emerge. First, the robustness issue of the MILP solutions needs to be addressed. Then, a natural extension of the presented work is to solve the same problem over longer horizons and to introduce a nonlinear battery model. Moreover, \(SoC_0\) could be optimized, with a terminal cost penalizing deviations from a final reference. Lastly, developing specialized interior point solvers exploiting the resulting bilinear structure could improve the efficiency of the NLP formulation. These extensions would enable a more meaningful comparison between the NLP and MILP formulations.
\bibliography{ifacconf}
\end{document}